\documentclass[a4paper,english]{article}

\usepackage{trcover}
\usepackage{psfrag,graphicx,amsfonts,amsmath,amssymb,amstext,amsthm,url,pdfsync,a4wide,enumerate,hyperref,multirow,xcolor}
\newcommand{\R}{\mathbb R}

\newcommand{\ps}[2]{\langle {#1} , {#2} \rangle }
\newcommand{\eps}{\varepsilon}

\newtheorem*{thank*}{Acknowledgments}
\newtheorem{definition}{Definition}[section]
\newtheorem*{definition*}{Definition}

\newtheorem*{example*}{Example}
\newtheorem*{application*}{Application}

\newtheorem*{lemma*}{Lemma}

\newtheorem*{note*}{Note}
\newtheorem*{notation*}{Notation}
\newtheorem{proposition}{Proposition}[section]
\newtheorem*{proposition*}{Proposition}
\newtheorem{remark}{Remark}[section]
\newtheorem*{remark*}{Remark}
\newtheorem{theorem}{Theorem}[section]
\newtheorem*{theorem*}{Theorem}
\newtheorem{corollary}{Corollary}[section]

\title{Solving non-monotone equilibrium problems via a DIRECT-type approach}

\trnum{}

\author{Stefano Lucidi\thanks{Department of Computer, Control and Management Engineering ``Antonio Ruberti'', Sapienza University of Rome, via Ariosto, 25, 00185 Roma, Italy, email: \texttt{lucidi@diag.uniroma1.it}.
	}
	\and
	Mauro Passacantando\thanks{Department of Computer Science, University of Pisa, 
		Largo B. Pontecorvo 3, 56127 Pisa, Italy, email: \texttt{mauro.passacantando@unipi.it}. 
	}
	\and
	Francesco Rinaldi\thanks{Department of Mathematics ``Tullio Levi-Civita'', University of Padova, via Trieste, 63, 35121 Padova, Italy, email: \texttt{rinaldi@math.unipd.it}. 
	}
}

\date{February 27, 2020}

\begin{document}

\maketitle

\begin{abstract}
A global optimization approach for solving non-monotone equilibrium problems (EPs) is proposed. The class of (regularized) gap functions is used to reformulate any EP as a constrained global optimization program and some  bounds on the Lipschitz constant of such functions are provided. 
The proposed global optimization approach is a combination of an improved version of the \texttt{DIRECT} algorithm, which exploits local bounds of the Lipschitz constant of the objective function, with local minimizations. 
Unlike most existing solution methods for EPs, no monotonicity-type condition is assumed in this paper. Preliminary numerical results on several classes of EPs show the effectiveness of the approach.
\end{abstract}

\

\noindent\textbf{Keywords:} Equilibrium problem, Gap function, Global optimization, \texttt{DIRECT} algorithm.

\


\section{Introduction}

Given a bifunction $f: \R^n  \times \R^n \to \R$ and a closed convex set $\mathcal{C} \subseteq \R^n$, we consider the following equilibrium problem:
\begin{equation} \tag{EP}
\text{find $x^* \in \mathcal{C}$ such that $f(x^*,y) \geq 0$, for any $y \in \mathcal{C}$.}
\end{equation}
This framework is a general mathematical model which includes several problems such as scalar and vector optimization, variational inequality (VI), complementarity, saddle point, Nash equilibrium problems in noncooperative games and inverse optimization~\cite{bigi2013existence,blum1994optimization}.

Several classes of iterative methods to solve EPs have been proposed in the literature: fixed point approaches~\cite{mastroeni2003auxiliary,muu2009regularization,nguyen2009bundle}, 
extragradient methods~\cite{flaam1996equilibrium,langenberg2013interior1,langenberg2013interior2,quoc2008extragradient}, 
descent algorithms~\cite{chadli2004descent,dilorenzo2014convergent,konnov2003Dgap,mastroeni2003gap}, 
proximal point methods~\cite{burachik2012generalized,konnov2003application,moudafi1999proximal}.
All these approaches need, directly or indirectly, some monotonicity-type assumption on the bifunction $f$ (e.g. strong or weak monotonicity, pseudo\-monotonicity, $\nabla$-monotonicity, etc.) in order to guarantee the convergence to a solution of (EP). 
On the other hand, it is well known that, without any need of monotonicity-type assumptions on $f$, (EP) can be reformulated as an equivalent global optimization problem via merit functions~\cite{pappalardo2016merit}. 
This fact suggests to use global optimization approaches to solve non-monotone EPs. 
Global optimization techniques have been considered in the literature only for two special classes of EPs:  linear complementarity problems~\cite{al1987implicit,pardalos1988global} and VI problems (a branch and bound method was proposed in~\cite{majig2009global} and a meta-heuristic algorithm in~\cite{majig2007hybrid}).

In this paper, we propose a \texttt{DIRECT}-type global optimization approach for solving general EPs, without assuming any monotonicity-type condition on $f$. 
In particular, we first reformulate (EP) as a global optimization problem via the well-known gap functions. 
We analyze the Lipschitz continuity of gap functions and give simple estimates of the Lipschitz constant for some special classes of EPs. 
Then, we combine the improved version of the \texttt{DIRECT} algorithm developed in~\cite{di2016direct}, which exploits local bounds of the Lipschitz constant of the objective function, with local searches to find a global minimum point of the gap function, i.e., a solution of (EP).
Finally, we show the effectiveness of our approach with some preliminary numerical experiments on instances coming from the literature and randomly generated instances.

The rest of the paper is organized as follows. In Section~\ref{s:background}, we recall the definition and the main properties of gap functions for (EP).
In Section~\ref{s:Lipschitz-gap}, we provide some general results on the Lipschitz continuity of gap functions and give explicit bounds of the Lipschitz constant for three classes of problems: affine VIs, VIs with trigonometric terms and affine EPs. 
Section~\ref{s:direct} presents the \texttt{DIRECT}-type global optimization approach and recalls the convergence properties of both the standard version of the \texttt{DIRECT} algorithm and its improved version proposed in~\cite{di2016direct}.
Section~\ref{s:results} reports the results of some preliminary numerical tests and shows that the improved version of \texttt{DIRECT} is more efficient than its standard version on most of the considered instances. 
Conclusions are finally drawn in Section~\ref{s:conclusions}. 

Throughout the paper we will assume that the feasible set $\mathcal{C}$ is compact, the bifunction $f$ is continuous, $f(x,\cdot)$ is convex and $f(x,x)=0$ for any $x \in \mathcal{C}$. It is well known that under these assumptions the existence of at least one solution of (EP) is guaranteed (see, e.g.,\cite{fan1972minimax}).  


\section{Preliminary background}
\label{s:background}

Merit functions allow reformulating (EP) as a global optimization problem, whose optimal value is known a priori. Several classes of merit functions for EPs have been introduced in the literature in the last two decades~\cite{pappalardo2016merit}. In this paper, we focus on the class of gap functions.

\begin{theorem}\label{t:gap}\cite{mastroeni2003gap}
	For any $\alpha \geq 0$ the gap function
	\begin{equation}\label{e:gap}
	\varphi_\alpha(x) := \max_{y \in \mathcal{C}} \left[ -f(x,y) - \frac{\alpha}{2}\,\|y-x\|^2 \right]
	\end{equation}
	has the following properties:
	\begin{enumerate}[a)]
		\item $\varphi_\alpha(x) \geq 0$ for any $x \in \mathcal{C}$;
		\item $x^*$ solves (EP) if and only if $x^* \in \mathcal{C}$ and $\varphi_\alpha(x^*)=0$;
		\item If $\alpha >0$ and $f$ is continuously differentiable on $\R^n \times \R^n$, then $\varphi_\alpha$ is continuously differentiable on $\R^n$ with 
		\begin{equation}\label{e:gap-grad}
		\nabla \varphi_\alpha(x) = -\nabla_1 f(x,y_\alpha(x)) - \alpha (x - y_\alpha(x)),
		\end{equation}
		where $\nabla_1 f(x,y)$ denotes the gradient of $f(\cdot,y)$ at $x$ and $y_\alpha(x)$ is the unique maximizer of problem in~\eqref{e:gap}.
	\end{enumerate}
\end{theorem}
Therefore, the solutions of (EP) coincide with the global minimum points of the optimization problem
\begin{equation}\label{e:minprob}
\begin{array}{rl}
\min & \varphi_\alpha(x)
\\[1mm]
& x \in \mathcal{C},
\end{array}
\end{equation} 
whose global minimum value is zero.
We remark that evaluating the gap function $\varphi_\alpha$ at some point $x$ consists in maximizing a concave (when $\alpha=0$) or strongly concave (when $\alpha>0$) function over the set $\mathcal{C}$. 
Moreover, the regularization term $\|y-x\|^2$ can be replaced by a more general bifunction satisfying suitable conditions (see~\cite{mastroeni2003gap}).

Several descent methods based on the gap function $\varphi_\alpha$ have been developed in the literature for solving EPs (see, e.g.,\cite{chadli2004descent,dilorenzo2014convergent,mastroeni2003gap}. However, their convergence to a solution of (EP) is guaranteed provided that some monotonicity-type assumption on the bifunction $f$ is assumed. 
In this paper, we propose a global optimization approach for solving problem~\eqref{e:minprob} that is not based on any monotonicity-type condition on $f$. More specifically, we consider a \texttt{DIRECT}-type method (see, e.g., \cite{jones1993lipschitzian}) with local searches. \texttt{DIRECT} (\texttt{DI}vide \texttt{RECT}angle) is a partitioning  strategy that  samples  points  in  the  domain  and  uses only objective function evaluations to  decide  what  to  do  next.  The boosted version we use here, called $\bar L$-\texttt{DIRECT} and first proposed in~\cite{di2016direct}, exploits overestimates of the Lipschitz constant related to the objective function
to improve the way the subsets to be further partitioned are selected. As we will see in the next section, this choice is well-suited to our problem.  Indeed, when our problem has some specific structure, an overestimate of the Lipschitz constant for the function $\varphi_\alpha$ can be easily calculated.

In the rest of the paper, we will consider the class of EPs where the bifunction
$$
f(x,y)=\ps{F(x,y)}{y-x}
$$ 
for some map $F: \R^n  \times \R^n \to \R^n$. This class of EPs includes two important particular cases: (i) VIs, where the map $F$ only depends on the variable $x$ and (ii) affine EPs, where $F(x,y)=Px+Qy+r$ for some  $P,Q \in \R^{n \times n}$ and $r \in \R^n$. Notice that Nash EPs in noncooperative games with quadratic cost functions are an interesting particular
case of affine EPs (see, e.g.,\cite{bigi2015twelve}).


\section{Lipschitz continuity of gap functions}
\label{s:Lipschitz-gap}

In this section, we provide some general results on the Lipschitz continuity of gap function $\varphi_\alpha$ and show some simple estimates of its Lipschitz constant for three special classes of EPs.
The knowledge of the Lipschitz constant of $\varphi_\alpha$ will be exploited by the global optimization approach described in Section~\ref{s:direct} for solving problem~\eqref{e:minprob}.

\begin{theorem}\label{t:Lip1}
	Suppose that $B \subseteq \R^n$ is compact, $F$ is continuous on $\R^n \times \R^n$ and $F(\cdot,y$) is Lipschitz continuous on $B$, uniformly with respect to $y$, with constant $L_F$. 
	Then, for any $\alpha \geq 0$ the function $\varphi_\alpha$ is Lipschitz continuous on $B$ with constant 
	$$
	L_1 + L_2\,L_F + \alpha\,L_2,
	$$
	where
	\begin{equation} \label{eq:L1eL2}
	L_1 = \max_{x \in B,\, y \in \mathcal{C}} \|F(x,y)\|,
	\qquad
	L_2 = \max_{x \in B,\, y \in \mathcal{C}} \|x - y\|.
	\end{equation}	
\end{theorem}

\begin{proof}
	If $x,y \in B$, then the following chain of equalities and inequalities holds:
	\begin{align*}
	\varphi_\alpha(x)-\varphi_\alpha(y) 
	& = \max_{z \in \mathcal{C}} \left[ \ps{F(x,z)}{x-z} - \frac{\alpha}{2}\|x-z\|^2 \right]
	- \max_{z \in \mathcal{C}} \left[ \ps{F(y,z)}{y-z} - \frac{\alpha}{2}\|y-z\|^2 \right] \\
	& \le \max_{z \in \mathcal{C}} \left[ \ps{F(x,z)}{x-z} - \ps{F(y,z)}{y-z} - \frac{\alpha}{2}\|x-z\|^2 + \frac{\alpha}{2}\|y-z\|^2 \right] \\
	& = \max_{z \in \mathcal{C}} \left[ \ps{F(x,z)-F(y,z)}{x-z} + \ps{F(y,z)}{x-y} + \frac{\alpha}{2}\ps{y-x}{y-z+x-z} \right] \\
	& \leq \max_{z \in \mathcal{C}} \left[ \|F(x,z)-F(y,z)\|\|x-z\| + \|F(y,z)\|\|x-y\| \right.
	\\
	& \qquad \left. + \frac{\alpha}{2}\|y-x\|(\|y-z\|+\|x-z\|) \right] 
	\\
	& \leq L_F\|x-y\|\,(\max_{z \in \mathcal{C}} \|x-z\|) + L_1\|x-y\| 
	\\
	& \qquad 
	+ \frac{\alpha}{2}\|y-x\| \left[ \max_{z \in \mathcal{C}} \|y-z\| + \max_{z \in \mathcal{C}} \|x-z\| \right] 
	\\
	& \le (L_1 + L_2\,L_F + \alpha\,L_2)\,\|x-y\|,
	\end{align*}
	where the second inequality follows from the Cauchy-Schwarz inequality, the third one from the Lipschitz continuity of $F$ and the last one from the definition of $L_2$. 
\end{proof}

\begin{remark}
	\autoref{t:Lip1} is a generalization of Lemma 2.1 proved in~\cite{majig2009global}, which provides an estimate of the Lipschitz constant of the gap function $\varphi_0$ for a VI with Lipschitz continuous operator. 
	In fact, when (EP) reduces to a VI, the regularization parameter $\alpha=0$ and the set $B=\mathcal{C}$, then the value of the Lipschitz constant given in \autoref{t:Lip1} coincides with that given in~\cite[Lemma 2.1]{majig2009global}.
\end{remark}

A further estimate of the Lipschitz constant of $\varphi_\alpha$, with $\alpha>0$, can be obtained provided that the map $F$ is smooth.

\begin{theorem}\label{t:Lip2}
	Suppose that $B \subseteq \R^n$ is a convex compact set and $F$ is continuously differentiable on $\R^n \times \R^n$. 
	Then, for any $\alpha > 0$ the function $\varphi_\alpha$ is Lipschitz continuous on $B$ with constant 
	$$
	L_1 + L_2\,L_3(\alpha),
	$$
	where $L_1$ and $L_2$ are defined in~\eqref{eq:L1eL2} and
	\begin{equation} \label{e:L3}
	L_3(\alpha) = 
	\max\limits_{x \in B,\, y \in \mathcal{C}} \| \alpha\,I - \nabla_1 F(x,y) \|,
	\end{equation}
	where $\nabla_1 F(x,y)$ denotes the Jacobian matrix of $F(\cdot,y)$ at $x$.
\end{theorem} 

\begin{proof}
	\autoref{t:gap} guarantees that $\varphi_\alpha$ is continuously differentiable on $\R^n$ with
	$$
	\nabla \varphi_\alpha(x) 
	= 
	F(x,y_\alpha(x)) +
	[ \alpha\,I - \nabla_1 F(x,y_\alpha(x))^T ] ( y_\alpha(x) - x ),
	\qquad x \in \R^n,
	$$
	where
	$$
	y_\alpha(x) 
	= 
	\arg\max_{y \in \mathcal{C}}
	\left[ \ps{F(x,y)}{x-y} - \frac{\alpha}{2}\,\|y-x\|^2 \right].
	$$
	Let $u, v \in B$. The mean value theorem guarantees that there exists $\xi \in (0,1)$ such that
	$$
	\varphi_\alpha(u) - \varphi_\alpha(v) = \ps{\nabla \varphi_\alpha(z)}{u-v},
	$$
	where $z:=\xi u + (1-\xi) v \in B$. Therefore, we get
	$$
	\begin{array}{rl}
	| \varphi_\alpha(u) - \varphi_\alpha(v) | 
	& \leq \| \nabla \varphi_\alpha(z) \| \, \|u-v\| 
	\\[2mm]
	& \leq \left[ \|F(z,y_\alpha(z))\| + \| \alpha\,I - \nabla_1 F(z,y_\alpha(z))^T \| \, \| y_\alpha(z) - z \| \right] \, \|u-v\| 
	\\[2mm]
	& \leq [ L_1 + L_2\,L_3(\alpha) ] \, \|u-v\|. 
	\end{array}
	$$
\end{proof}

In the special case of a VI defined by a smooth map, a third estimate of the Lipschitz constant of $\varphi_\alpha$ can be proved.

\begin{theorem}\label{t:Lip3}
	Suppose that (EP) is a VI, i.e., $f(x,y)=\ps{F(x)}{y-x}$ for some continuously differentiable map $F:\R^n \to \R^n$. 
	If $B \subseteq \R^n$ is a convex compact set such that $B \subseteq \mathcal{C}$, then for any $\alpha >0$ the function $\varphi_\alpha$ is Lipschitz continuous on $B$ with constant 
	$$
	L_1 + \alpha^{-1}\,L_1\,L_3(\alpha),
	$$
	where $L_1$ and $L_3(\alpha)$, defined in~\eqref{eq:L1eL2} and \eqref{e:L3} respectively, in this special case are equal to
	\begin{equation*}
	L_1 = \max_{x \in B} \| F(x) \|,
	\qquad
	L_3(\alpha) = \max\limits_{x \in B} \| \alpha\,I - \nabla F(x) \|.
	\end{equation*}
\end{theorem} 

\begin{proof}
	\autoref{t:gap} guarantees that $\varphi_\alpha$ is continuously differentiable and
	$$
	\nabla \varphi_\alpha(x) 
	= 
	F(x) + [ \alpha\,I - \nabla F(x)^T] ( y_\alpha(x) - x ),
	\qquad x \in \R^n,
	$$
	with
	$$
	y_\alpha(x) 
	= 
	P_\mathcal{C}(x - \alpha^{-1} F(x)),
	$$
	where $P_\mathcal{C}$ denotes the Euclidean projection on the set $\mathcal{C}$. 
	If $u, v \in B$, then the mean value theorem implies 
	$$
	\varphi_\alpha(u) - \varphi_\alpha(v) = \ps{\nabla \varphi_\alpha(z)}{u-v},
	$$
	where $z:=\xi u + (1-\xi) v$ for some $\xi \in (0,1)$. Therefore, we get
	$$
	\begin{array}{rl}
	| \varphi_\alpha(u) - \varphi_\alpha(v) | 
	& \leq \| \nabla \varphi_\alpha(z) \| \, \|u-v\| 
	\\[2mm]
	& \leq \left[ \|F(z)\| + \| \alpha\,I - \nabla F(z)^T \| \, \| y_\alpha(z) - z \| \right] \, \|u-v\| 
	\\[2mm]
	& = \left[ \|F(z)\| + \| \alpha\,I - \nabla F(z) \| \, \| P_\mathcal{C}(z-\alpha^{-1}F(z)) - P_\mathcal{C}(z) \| \right] \, \|u-v\| 
	\\[2mm]
	& \leq \left[ \|F(z)\| + \| \alpha\,I - \nabla F(z) \| \, \| z-\alpha^{-1}F(z) - z \| \right] \, \|u-v\| 
	\\[2mm]
	& \leq [L_1 + \alpha^{-1}\,L_1\,L_3(\alpha)] \, \|u-v\|,
	\end{array}
	$$
	where the third inequality holds since the projection map $P_\mathcal{C}$ is nonexpansive, i.e., $\|P_\mathcal{C}(x)-P_\mathcal{C}(y)\| \leq \|x-y\|$ holds for any $x, y \in \R^n$. 
\end{proof}

In the rest of this section we analyze the Lipschitz constant of $\varphi_\alpha$ for some special classes of EPs.


\subsection{Affine VIs defined on a box} 

Suppose that (EP) is a VI defined by an affine operator $F(x)=Px+r$, for some $P \in \R^{n \times n}$ and $r \in \R^n$, over a box $\mathcal{C}=[l,u]$, where $l, u \in \R^n$. 
Consider a box $B=[a,b]$, where $a,b \in \R^n$, such that $B \subseteq \mathcal{C}$, i.e., $l \le a \le b \le u$. 
Then, \autoref{t:Lip1}, \autoref{t:Lip2} and \autoref{t:Lip3} guarantee that $\varphi_0$ is Lipschitz continuous on $B$ with constant
\begin{equation}\label{e:VIaff_phi0}
L_1 + L_2\,L_F,
\end{equation}
while, for any $\alpha >0$, $\varphi_\alpha$ is Lipschitz continuous on $B$ with constant
\begin{equation}\label{e:VIaff_phialfa}
\min \left\{ 
L_1+L_2\,L_F+\alpha\,L_2,
\quad 
L_1 + L_2\,L_3(\alpha),
\quad
L_1 + \alpha^{-1}\,L_1\,L_3(\alpha) 
\right\}.
\end{equation}
We now show that the exact values (or upper bound) of the constants involved in the above formulas can be easily computed.

\

\noindent\textbf{Estimate of $L_1$.} The exact value of $L_1$ is
$$
L_1 = \max_{x \in B} \|Px+r\| = \max_{x \in vert(B)} \|Px+r\|,
$$
where $vert(B)$ denotes the set of vertices of $B$. Such a evaluation can be computationally expensive since the vertices of $B$ are exponentially many with respect to the number of variables. 
However, the following upper bounds for $L_1$ can be easily computed. If we denote by $P^+$ the Moore-Penrose pseudoinverse matrix of $P$, then we get
\begin{align*}
L_1 & = \max_{x \in B} \| Px+r \| 
\\
& = \max_{a \le x \le b} \|P(x+P^+ r) + (I-PP^+)r\| 
\\
& \leq \|(I-PP^+)r\| + \max_{a \le x \le b} \|P(x+P^+ r)\| 
\\
& \leq \|(I-PP^+)r\| + \|P\|\,\max_{a \le x \le b} \|x+P^+ r\| 
\\
& = \|(I-PP^+)r\| + \|P\|\,\|c(a,b)\| 
\\
& := L_1',
\end{align*}
where the $i$-th component of the vector $c(a,b)$ is defined as $c_i(a,b) = \max\{ |(a+P^+ r)_i| ,\ |(b+P^+ r)_i| \}$. Moreover, the following simple upper bounds hold:
$$
L_1 = \max_{a \le x \le b} \|P(x-a) + Pa+r\| 
\leq 
\|Pa+r\| + \|P\|\,\|b-a\|:=L_1^{''},
$$
$$
L_1 = \max_{a \le x \le b} \|P(x-b) + Pb+r\| 
\leq 
\|Pb+r\| + \|P\|\,\|b-a\|:=L_1^{'''}.
$$
Therefore, we have 
\begin{equation}\label{e:L1tilde}
L_1 \leq \widetilde{L}_1(P,r,a,b),
\qquad
\text{where}
\qquad
\widetilde{L}_1(P,r,a,b):=\min\{ L_1', L_1^{''}, L_1^{'''} \}.
\end{equation}

\begin{remark}
	In~\cite{majig2009global} the following upper bound for $L_1$ is given:
	$$
	L_1 \leq  \|P\|\,\|c(a,b)\|. 
	$$
	We remark that this inequality is not true in general, as the following counterexample shows. Let $n=2$,
	$$
	P=\left(\begin{array}{ccc} 1 & & 1 \\ 0 & & 0 \end{array} \right),
	\qquad
	r=\left(\begin{array}{c} 0 \\ v \end{array} \right),
	\ \text{ with $v \neq 0$, }
	\qquad
	a=\left(\begin{array}{c} 0 \\ 0 \end{array} \right),
	\qquad
	b=\left(\begin{array}{c} 1 \\ 1 \end{array} \right).
	$$
	Then, it is easy to check that $\|P\|=\sqrt{2}$ holds and the pseudoinverse of $P$ is
	$$
	P^+ = \left(\begin{array}{ccc} 1/2 & & 0 \\ 1/2 & & 0 \end{array} \right),
	$$
	hence $c_i(a,b)=\max\{|a_i|, \, |b_i|\}=1$ for $i=1,2$. Therefore, 
	$\|P\| \|c(a,b)\| = 2$. On the other hand, 
	$$
	L_1 
	= \max_{x \in B} \|Px+r\|
	= \max_{0 \le x \le 1} \|(x_1+x_2,v)\|
	= \|(2,v)\|
	= \sqrt{4+v^2} 
	> 2
	= \|P\| \|c(a,b)\|.
	$$
\end{remark}

\

\noindent\textbf{Estimate of $L_2$.}  The exact value of $L_2 = \max_{x \in B,\, y \in \mathcal{C}} \|x - y\|$ can be computed by solving $n$ independent optimization problems of the form
$$
\mathop{\max_{a_i \leq x_i \leq b_i}}_{l_i \leq y_i \leq u_i} (x_i-y_i)^2
= \max\{ (u_i-a_i)^2 , (l_i-b_i)^2 \},
$$
for $i=1,\dots,n$. Therefore, we have 
\begin{equation}\label{e:L2}
L_2 = \sqrt{\sum_{i=1}^n \max\{ (u_i-a_i)^2 , (l_i-b_i)^2 \}} .
\end{equation}

\

\noindent\textbf{Estimates of $L_3(\alpha)$ and $L_F$.} It is easy to check that $L_3(\alpha) = \| \alpha\,I - P \|$ and $L_F = \| P \|$.


\subsection{VIs with Trigonometric terms defined on a box} 

Suppose that (EP) is a VI defined over a box $\mathcal{C}=[l,u]$, with an operator which is the sum of an affine map and a trigonometric map, i.e.,
$$
F(x) = Px + r + T(x),
$$
where $T_i(x) = w_i \sin(v_i x_i)$, for $i=1,\dots,n$, $P \in \R^{n \times n}$ and $r, v, w \in \R^n$ with $v, w >0$. 
Consider a box $B=[a,b] \subseteq \mathcal{C}$, i.e., $l \leq a \leq b \leq u$. 
Then, \autoref{t:Lip1}, \autoref{t:Lip2} and \autoref{t:Lip3} imply that $\varphi_0$ is Lipschitz continuous on $B$ with constant~\eqref{e:VIaff_phi0},
while $\varphi_\alpha$, for any $\alpha >0$, is Lipschitz continuous on $B$ with constant~\eqref{e:VIaff_phialfa}.

\

\noindent\textbf{Estimate of $L_1$.} An upper bound for $L_1$ can be computed as follows:
$$
\begin{array}{rl}
L_1 & = \max\limits_{x \in B} \| Px + r + T(x) \| 
\\[3mm]
& \leq \max\limits_{x \in B} \|Px+r\| + \max\limits_{x \in B} \| T(x) \| 
\\[3mm]
& \leq \widetilde{L}_1(P,r,a,b) + \| w \|,
\end{array}
$$
where $\widetilde{L}_1(P,r,a,b)$ is defined in~\eqref{e:L1tilde}.

\

\noindent\textbf{Estimate of $L_2$.} Since $L_2$ only depends on the $B$ and $\mathcal{C}$, its exact value is given by~\eqref{e:L2}.

\

\noindent\textbf{Estimate of $L_3(\alpha)$.} The following upper bound can be obtained: 
$$
L_3(\alpha)
= 
\max_{x \in B} \| \alpha\,I - P - \nabla T(x) \|
\le
\|\alpha\,I - P\| + \max_{x \in B} \|\nabla T(x)\|.  
$$
The Jacobian matrix $\nabla T(x)$ is diagonal with
$$
[\nabla T(x)]_{ii} = w_i v_i \cos(v_i x_i),
\qquad i = 1,\dots,n,
$$
hence, for any $x \in B$ we get
$$
\|\nabla T(x)\|
=
\sqrt{\left[ \lambda_{max}(\nabla T(x)) \right]^2}
=
\sqrt{\left[ \max_{1 \leq i \leq n} \{ w_i v_i \cos(v_i x_i) \} \right]^2}
\leq
\sqrt{\left[ \max_{1 \leq i \leq n} \{ w_i v_i \} \right]^2}
=
\max_{1 \leq i \leq n} \{ w_i v_i \},
$$
where $\lambda_{max}(\nabla T(x))$ denotes the maximum eigenvalue of $\nabla T(x)$. Therefore, we have
$$
L_3(\alpha) \leq \|\alpha\,I - P\| + \max_{1 \leq i \leq n} \{ w_i v_i \}.
$$

\

\noindent\textbf{Estimate of $L_F$.} The Lipschitz constant of $F$ can be estimated as follows:
$$
\begin{array}{rl}
\|F(x)-F(z)\| 
&   = \|P(x-z) + T(x)-T(z)\| \\
& \le \|P\| \|x-z\| + \|T(x)-T(z)\| \\
&   = \|P\| \|x-z\| + \sqrt{\sum\limits_{i=1}^n w_i^2 
	[ \sin(v_i x_i)-\sin(v_i z_i) ]^2} \\
& \le \|P\| \|x-z\| + \sqrt{\sum\limits_{i=1}^n w_i^2 v_i^2 (x_i-z_i)^2} \\
& \le \|P\| \|x-z\| + \sqrt{\left[ \max\limits_{1 \leq i \leq n} \{ w_i v_i \} \right]^2 \sum\limits_{i=1}^n (x_i-z_i)^2} 
\\
&   = \|P\| \|x-z\| + \max\limits_{1 \leq i \leq n} \{ w_i v_i \} \|x-z\| \\    
&   = ( \|P\| + \max\limits_{1 \leq i \leq n} \{ w_i v_i \} ) \|x-z\|,
\end{array}
$$
where the second inequality holds because the sine function is Lipschitz continuous with constant~1. Therefore, we have
$$
L_F \leq \|P\| + \max\limits_{1 \leq i \leq n} \{ w_i v_i \}.
$$


\subsection{Affine EPs defined on a box}

Suppose that (EP) is defined by an affine operator $F(x,y)=Px+Qy+r$, for some $P,Q \in \R^{n \times n}$ and $r \in \R^n$, over a box $\mathcal{C}=[l,u]$, where $l,u \in \R^n$. 
Consider a box $B=[a,b]$. 
Then, \autoref{t:Lip1} and \autoref{t:Lip2} imply that $\varphi_0$ is Lipschitz continuous on $B$ with constant~\eqref{e:VIaff_phi0},
while, for any $\alpha >0$, $\varphi_\alpha$ is Lipschitz continuous on $B$ with constant
$$
\min \left\{ L_1+L_2\,L_F+\alpha\,L_2 ,\quad 
L_1 + L_2\,L_3 \right\}.
$$

\noindent\textbf{Estimate of $L_1$.} The following bound can be easily obtained: 
\begin{align*}
L_1 & = \max_{x \in B,\ y \in \mathcal{C}} \|Px+Qy+r\| \\
& \le \max_{x \in B,\ y \in \mathcal{C}} (\|Px\| + \|Qy+r\|) \\
& = \max_{x \in B} \|Px\| + \max_{y \in \mathcal{C}} \|Qy+r\| \\
& \leq \widetilde{L}_1(P,0,a,b) + \widetilde{L}_1(Q,r,l,u) \\
& := M_1. \\
\end{align*}
Similarly to the previous bound, we get
\begin{align*}
L_1 & = \max_{x \in B,\ y \in \mathcal{C}} \|Px+Qy+r\| \\
& \le \max_{x \in B,\ y \in \mathcal{C}} (\|Px+r\| + \|Qy\|) \\
& = \max_{x \in B} \|Px+r\| + \max_{y \in \mathcal{C}} \|Qy\| \\
& = \widetilde{L}_1(P,r,a,b) + \widetilde{L}_1(Q,0,l,u) \\
& := M_2. \\
\end{align*}
Finally, we have 
\begin{align*}
L_1 & = \max_{x \in B,\ y \in \mathcal{C}} \|Px+Qy+r\| \\
& \le \max_{x \in B,\ y \in \mathcal{C}} (\|Px+r/2\| + \|Qy+r/2\|) \\
& = \max_{x \in B} \|Px+r/2\| + \max_{y \in \mathcal{C}} \|Qy+r/2\| \\
& = \widetilde{L}_1(P,r/2,a,b) + \widetilde{L}_1(Q,r/2,l,u) \\
& := M_3, \\
\end{align*}
thus $L_1 \leq \min\{ M_1, M_2, M_3 \}$.

\

\noindent\textbf{Estimates of $L_2$, $L_3(\alpha)$ and $L_F$.} It is easy to check that $L_2$ is given by~\eqref{e:L2}, $L_3(\alpha) = \| \alpha\,I - P \|$ and $L_F = \| P \|$.


\section{The \texttt{\texttt{DIRECT}}-type  algorithms}
\label{s:direct}

We now describe a \texttt{DIRECT}-type approach to globally solve the optimization problem~\eqref{e:minprob} that is equivalent to (EP). In this section, we assume that the feasible region $\cal C$ is a box, i.e., \linebreak ${\cal C}=\{ x \in \R^n:\ l \leq x \leq u\}$.


More specifically, we focus on partition based algorithms, a class of methods with  both interesting theoretical properties and  efficient computational behavior,
and explain why those algorithms represent a good option when dealing with  non-monotone EPs. We start by giving some useful details. 

Partition based methods produce a sequence of finer and finer partitions $\{\mathcal{H}_k\}$ of the feasible set~$\cal C$.
At each iteration $k$, the $k$-th
partition is described by:
$$
\mathcal{H}_k=\{{\cal C}^i: i\in I_k\},
$$
where
$$
{\cal C}^i=\{x\in \mathbb{R}^n: l^i\le x\le u^i\}, \quad l^i,u^i\in [l,u], \quad x^i=(l^i+u^i)/2.
$$
Then the next partition $\mathcal{H}_{k+1}$ is obtained by selecting and by further partitioning every element of  a
``particular'' subset  $\{{\cal C}^i: i\in I^*_k\}\subseteq\mathcal{H}_{k}$, where $I^*_k\subset I_k$.
A partition based algorithm is characterized by the rules used to generate the subset of indices $I^*_k$, and by the strategies applied to further partition the subsets~$\{{\cal C}^i: i\in I^*_k\}$.

In \cite{majig2009global},  the authors consider non-monotone VIs and use a Branch and Bound method  similar to the one described in \cite{horst2000introduction}
to tackle the considered global optimization problems.

Instead, as previously pointed out, we solve non-monotone EPs by means of an algorithm  derived from the well-known  \texttt{DIRECT} method (see, e.g., \cite{jones1993lipschitzian}). 
This approach, called $\bar L$-\texttt{DIRECT} and first proposed in~\cite{di2016direct}, differs from the standard version of \texttt{DIRECT} in the way the set of indices $I^*_k$ are defined.
In the  standard version of \texttt{DIRECT}, $I^*_k$ consists of the  indices related to the subsets satisfying the definition reported below:
\begin{definition}\label{potenz-ott-f}
	Given a partition $\mathcal{H}_k=\{ {\cal C}^i: i\in
	I_k\}$  of ${\cal C}$ and  a scalar $\eps>0$, a  subset ${\cal C}^h$ is  \emph{potentially optimal with respect to the function $\varphi_\alpha$} if a constant
	$\bar L^h$ exists such that:
	\begin{eqnarray*}\label{rel:potopt1}
		&& \varphi_\alpha(x^h)-\frac{\bar L^h}{2}\|u^h-l^h \|\le \varphi_\alpha(x^i)-\frac{\bar
			L^h}{2}\|u^i-l^i \|,\qquad\qquad \forall\ i\in I_k,
		\\
		\label{rel:potopt21}&& \varphi_\alpha(x^h)-\frac{\bar L^h}{2}\|u^h-l^h \|\le \varphi_{\min}-\epsilon |\varphi_{\min}|,
	\end{eqnarray*}
	where
	\begin{equation}\label{fmin}\varphi_{\min}=\min_{i\in I_k}\varphi_\alpha\bigl(x^i\bigr).\end{equation}
\end{definition}
In the $\bar L$-\texttt{DIRECT} algorithm, $I^*_k$ is given by the indices related to those subsets satisfying:
\begin{definition}
	\label{def-strong-po}
	Given a partition $\mathcal{H}_k=\{ {\cal C}^i: i\in
	I_k\}$  of ${\cal C}$, a scalar $\eps>0$, a scalar $\eta
	>0$ and a scalar $\bar L>0$ , a subset  ${\cal C}^h$ is
	\emph{$\bar L$-potentially optimal with respect to the function $\varphi_\alpha$} if one of the following conditions is satisfied:
	\begin{itemize}
		\item[i)] A constant $\tilde L^h \in (0, \bar L)$ exists such that:
		\begin{eqnarray}
		&& \varphi_\alpha(x^h)-\frac{\tilde L^h}{2}\|u^h-l^h \|\le \varphi_\alpha(x^i)-\frac{\tilde
			L^h}{2}\|u^i-l^i \|,\qquad\qquad \forall\ i\in I_k,\\
		\label{rel:potopt22}&& \varphi_\alpha(x^h)-\frac{\tilde L^h}{2}\|u^h-l^h \|\le \varphi_{\min}-\epsilon \max\{|\varphi_{\min}|,  \eta\},
		\end{eqnarray}
		where $\varphi_{\min}$ is given by (\ref{fmin});
		\item[ii)] The following inequality holds:
		\begin{eqnarray}
		&& \varphi_\alpha(x^h)-\frac{\bar L}{2}\|u^h-l^h \|\le \varphi_\alpha(x^i)-\frac{\bar L}{2}\|u^i-l^i \|,\qquad\qquad \forall\ i\in I_k.
		\end{eqnarray}
	\end{itemize}
\end{definition}
\begin{remark}
	The difference between the two is that an  overestimate $\bar L$ of the Lipschitz constant is used in Definition \ref{def-strong-po}. This fact obviously enhances the way $\bar L$-\texttt{DIRECT}  selects the subsets to be partitioned.
\end{remark}

\begin{remark}
	As $\bar L\to\infty$, Definition \ref{def-strong-po} tends to Definition \ref{potenz-ott-f} and, hence, the strategy proposed in \cite{di2016direct} becomes the one proposed in \cite{jones1993lipschitzian}.
\end{remark}
We refer to  \cite{jones1993lipschitzian} and  \cite{di2016direct} for detailed descriptions and discussions of the \texttt{DIRECT} algorithm and the $\bar L$-\texttt{DIRECT} algorithm.
Similarly to  any partition-based method, the asymptotic behavior shown by the \texttt{DIRECT}  and the  $\bar L$-\texttt{DIRECT} algorithms is characterized by   the partition sequences they produce. Those   sequences can be  represented equivalently by infinite sequences  of nested subsets $\{{\cal C}^{i_k}\}$, defined as follows:
\begin{description}
	\item[] Given a set ${\cal C}^{i_k}$  at the iteration $k$, its  predecessor ${\cal C}^{i_{k-1}}$
	is the unique set belonging to the previous partition $\mathcal{H}_{k-1}=\{ {\cal C}^i: i\in
	I_{k-1}\}$ such that  ${\cal C}^{i_{k}}\subseteq {\cal C}^{i_{k-1}}$.
\end{description}
Then, the analysis of  theoretical properties of \texttt{DIRECT} algorithm and $\bar L$-\texttt{DIRECT} algorithm can be performed by studying the properties of the produced sequences $\{{\cal C}^{i_k}\}$.
The partitioning strategy used by the \texttt{DIRECT} algorithm and the $\bar L$-\texttt{DIRECT} algorithm guarantees (regardless of the particular choice of set $I_k^*$) that the produced sequences $\{{\cal C}^{i_k}\}$ satisfy one of the following properties (see \cite{DApuzzo}):
\begin{description}
	\item[-] {\cal Property 1}:\qquad an index $\bar k$ exists such that ${\cal C}^{i_{\bar k}}= {\cal C}^{i_{k}}$ for all $k\ge\bar k$;
	\item[-] {\cal Property 2}:\qquad $\displaystyle\bigcap_{k=0}^\infty {\cal C}^{i_{k}}=\{\bar x\},\quad \mbox{where}\quad \bar x\in {\cal C}$.
\end{description}
Then the so-called {\it everywhere dense convergence}  can be stated by the following proposition.
\begin{proposition} \label{denso-nonvinc} \texttt{DIRECT} algorithm has the following properties:
	\begin{itemize}
		\item[i)]  All the sequences of sets $\{{\cal C}^{i_k}\}$ produced  satisfy Property 2;
		\item[ii)]  For every $\tilde x
		\in {\cal C}$, the \texttt{DIRECT} algorithm produces a sequence of sets  $\{{\cal C}^{i_k}\}$ satisfying Property 2 and such that
		$$\bigcap_{k=0}^\infty {\cal C}^{i_{k}}=\{\tilde x\}.$$
	\end{itemize}
\end{proposition}
The properties of  the  $\bar L$-\texttt{DIRECT} algorithm also depend on the choice of the scalar $\bar L$  included in the definition of
$\bar L$-potentially optimal subsets. In particular, the following assumption can be introduced.

\noindent{\bf Assumption 1.} {\it For every global minimum point $ x^*$ of problem \eqref{e:minprob}, there exists an index $\bar k$ (possibly depending on $x^*$) such
	that, if \ ${\cal C}^{j_{\bar k}}\in \{ {\cal C}^i: i\in
	I_{\bar k}\}$ is the subset satisfying $x^*\in {\cal C}^{j_{\bar k}}$, then
	$$
	\bar L< L,
	$$
	where $L$ is the local Lipschitz constant of the function $\varphi_\alpha$ over the subset ${\cal C}^{j_{\bar k}}$.}

Now it possible to state the following result.
\begin{proposition}\label{prop-conv-VarDirect}
	If Assumption 1 holds, then $\bar L$-\texttt{DIRECT} algorithm has the following properties:
	\begin{description}
		\item{i)} Every  sequence of sets  $\{{\cal C}^{i_k}\}$ produced by the algorithm which satisfies Property 2 is such that
		$$
		\bigcap_{k=0}^\infty {\cal C}^{i_k}=\{x^*\},
		$$
		where $ x^*$ is a global minimum of problem \eqref{e:minprob};
		\item{ii)} For every global minimum  $ x^*$ of problem \eqref{e:minprob}, the algorithm produces a  sequence of sets  $\{{\cal C}^{i_k}\}$ satisfying Property 2 and
		$$
		\bigcap_{k=0}^\infty {\cal C}^{i_{k}}=\{ x^* \};
		$$
		\item{iii)} Let $\bar k$ be the index introduced in Assumption 1. Then, for all $k\ge \bar k$, the following inequality holds
		\begin{equation}\label{stopc}\varphi_\alpha(x^{h_k})- \varphi_\alpha^*\le  \frac{\bar L}{2}\|u^{h_k}-l^{h_k} \|,\end{equation}
		where the index ${h}_k$ is given by:
		$$ \varphi_\alpha(x^{h_k})-\frac{\bar L}{2}\|u^{h_k}-l^{h_k} \|= \min_{i\in I_k} \bigg\{\varphi_\alpha(x^i)-\frac{\bar L}{2} \|u^{i}-l^{i} \|\bigg\}.$$
	\end{description}
\end{proposition}
Points i) and ii) of the previous proposition guarantee that, as the number of iterations increases, $\bar L$-\texttt{DIRECT} generates points that are more and more clustered  around the global minima of problem~\eqref{e:minprob}. 
Point iii) gives a practical stopping criterion for the algorithm. The right-hand side of \eqref{stopc} indeed provides an optimality gap.

\begin{remark}
	Proposition \ref{prop-conv-VarDirect} highlights the main difference between the Branch and Bound algorithm used in \cite{majig2009global} and  the $\bar L$-\texttt{DIRECT} algorithm. 
	In order to guarantee convergence to a global minimum of the Branch and Bound,  an  overestimate for the Lipschitz constant of $\varphi_\alpha$ over the whole feasible set ${\cal C}$ is needed from the beginning. On the other hand,  convergence of the $\bar L$-\texttt{DIRECT} algorithm can be guaranteed by  an  overestimate of the local  Lipschitz constant of $\varphi_\alpha$ over the  subset ${\cal C}^{j_{\bar k}}$ (keep in mind that this local constant is usually much smaller than the global one). Furthermore,
	this overestimate  is needed only for sufficiently large values of the indices $k$.  Hence, the information obtained from the function values calculated in the first iterations of the algorithm can be exploited to get an overestimate of the required local Lipschitz constant.
\end{remark}

Proposition \ref{denso-nonvinc} and Proposition   \ref{prop-conv-VarDirect} imply the following corollary.

\begin{corollary}
	The \texttt{DIRECT} Algorithm and, if Assumption 1 holds, also the $\bar L$-\texttt{DIRECT} algorithm satisfy the following property:
	\begin{description}
		\item[] For every every global minimum point $ x^*$ of problem~\eqref{e:minprob} and for every neighborhood ${\cal B}(x^*)$ of $ x^*$, an index $\bar k$ exists such that both the algorithms produce a point $x^{i_{\bar k}}$ satisfying
		$$x^{i_{\bar k}}\in{\cal B}(x^*) .$$
	\end{description}
\end{corollary}

The previous result points out that the two \texttt{DIRECT}-based methods can be efficiently combined with local searches within a multistart strategy.


\section{Numerical Results}
\label{s:results}

In this section, we describe our numerical experience. The goal is twofold: on the one side, we would like to see how \texttt{DIRECT} strategies behave on this class of problems; on the other side, we would like to understand the importance of embedding the Lipschitz constant estimates in those algorithmic schemes. 
We thus consider two different algorithms in the experiments:
\begin{itemize}
	\item[$\bullet$] \texttt{DIRECT}: the standard version of the method with local searches;
	\item[$\bullet$] \texttt{$\bar L$-DIRECT}: the modified version with Lipschitz constant estimates and local searches.
\end{itemize}
In both cases we used the SDBOX algorithm~\cite{lucidi2002derivative} to perform the local search. All 
algorithms were implemented in Matlab and tests were performed with Matlab v2019b.
We first considered randomly generated instances for two different classes of problems, that is affine VIs 
and VIs with trigonometric terms. In the analysis of those randomly generated instances we used performance
and data profiles~\cite{more2009benchmarking} with a gate parameter $\tau=10^{-3}$.
Then, we considered 5 affine VI problems coming from the literature. In all the experiments, we considered the gap function $\varphi_\alpha$ defined in~\eqref{e:gap} with $\alpha=1$.
The detailed results are reported in the next subsections.


\subsection{Results on Randomly Generated Affine VIs} 

\begin{figure}[htbp]
	\begin{center}
		\psfrag{Direct}{\scalebox{.45}{\texttt{DIRECT}}}
		\psfrag{Direct-lip}{\scalebox{.45}{$\bar L$-\texttt{DIRECT}}}
		\includegraphics[width=0.69\textwidth]{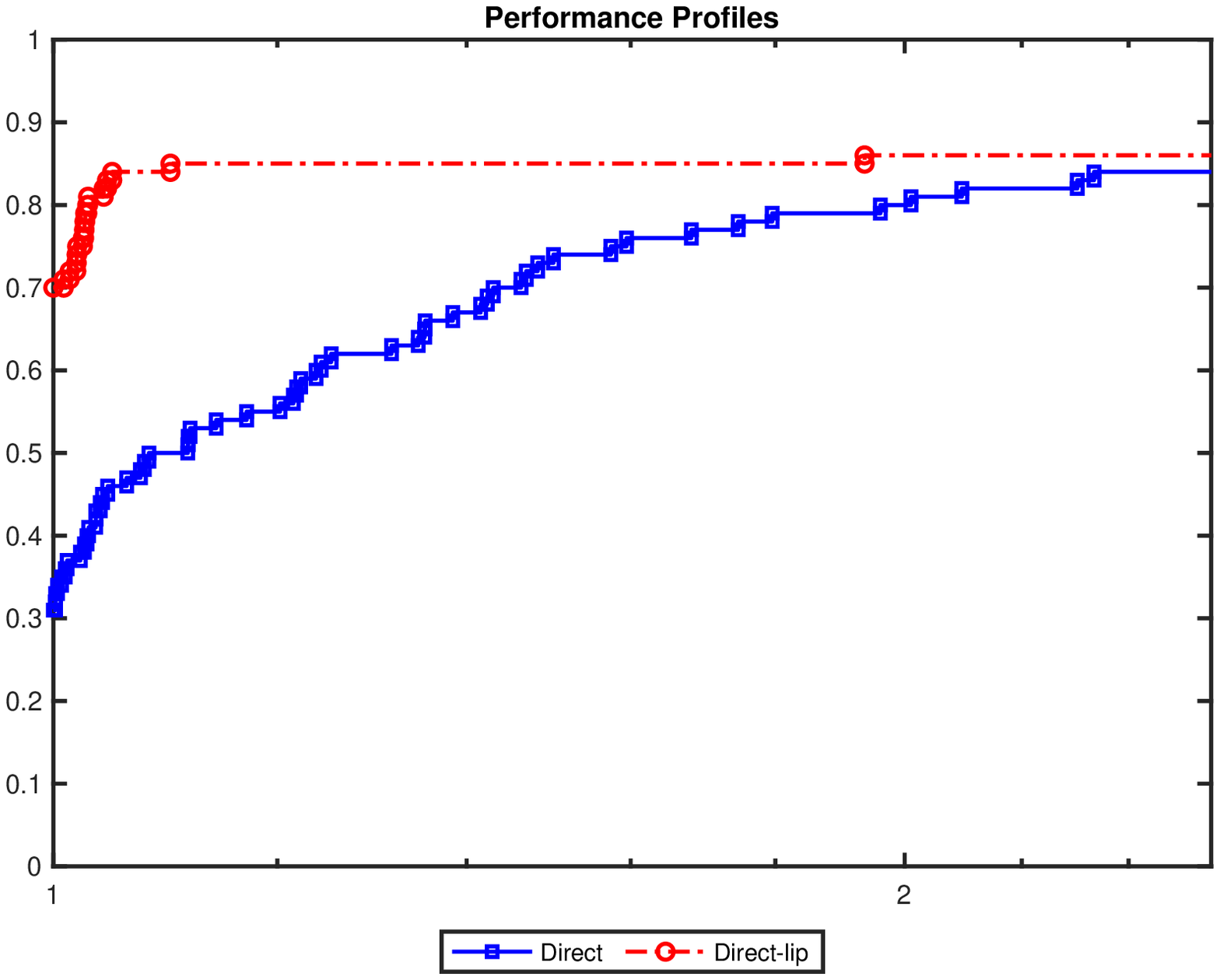}
		\\
		\includegraphics[width=0.69\textwidth]{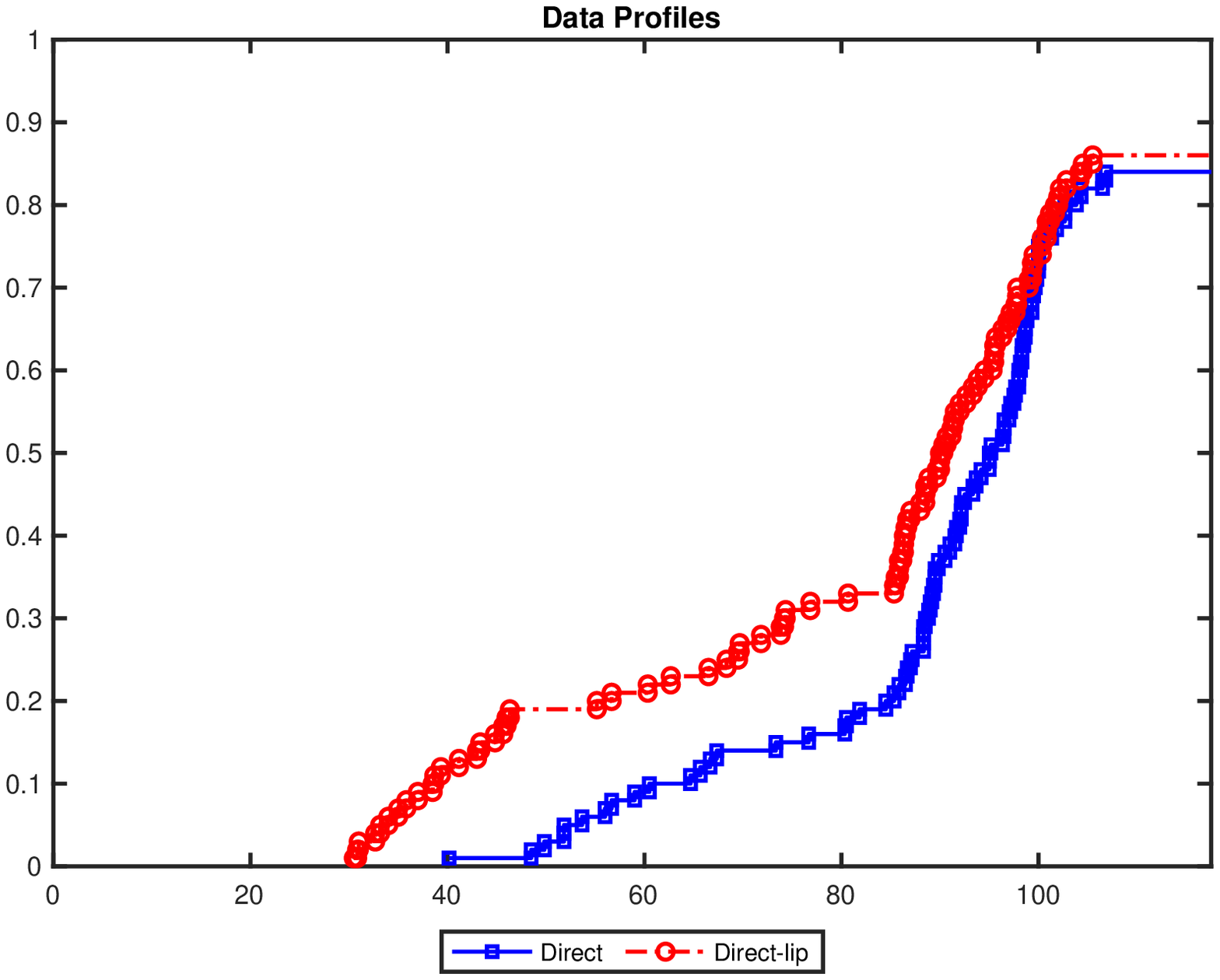}
	\end{center}
	\caption{Performance and data profiles for randomly generated affine VI problems.}
	\label{fig:2}       
\end{figure}

We now describe in depth the results obtained on randomly generated affine VI problems. We generated 100 instances with 5 variables.
For each instance, the affine operator $F(x)=Px+r$ was randomly built by choosing a matrix $P$ with uniformly distributed random
numbers in the interval $[0,3]$ and a vector $r$  with uniformly distributed random
numbers in the interval $[-2,2]$. The box constraints $\{x \in \R^n:\ l\leq x\leq u\}$ were generated by 
considering two vectors $l$ and $u$  with uniformly distributed random numbers in the interval $[-2,0]$ and $[1,3]$, respectively.
We gave a budget of 600 function evaluations  to the considered algorithms (500 for the \texttt{DIRECT} strategies and 100 for the local search).
Performance and data profiles are reported in~\autoref{fig:2}. The  performance profile plot shows that
the \texttt{$\bar L$-DIRECT} (red line) is both much more efficient than \texttt{DIRECT} (blue line), since it gives better performance and satisfies the stopping condition 
with a smaller number of function evaluations for the $70\%$ of the instances, and  more reliable (indeed, the percentage of problems that can be solved with the available budget of function evaluations is higher). If we observe the data profiles, we can further see that \texttt{$\bar L$-DIRECT} solves a higher percentage of problems no matter what the budget used is. 

\subsection{Results on Randomly Generated VIs with Trigonometric Terms}

\begin{figure}[htbp]
	\begin{center}
		\psfrag{Direct}{\scalebox{.45}{\texttt{DIRECT}}}
		\psfrag{Direct-lip}{\scalebox{.45}{$\bar L$-\texttt{DIRECT}}}
		\includegraphics[width=0.69\textwidth]{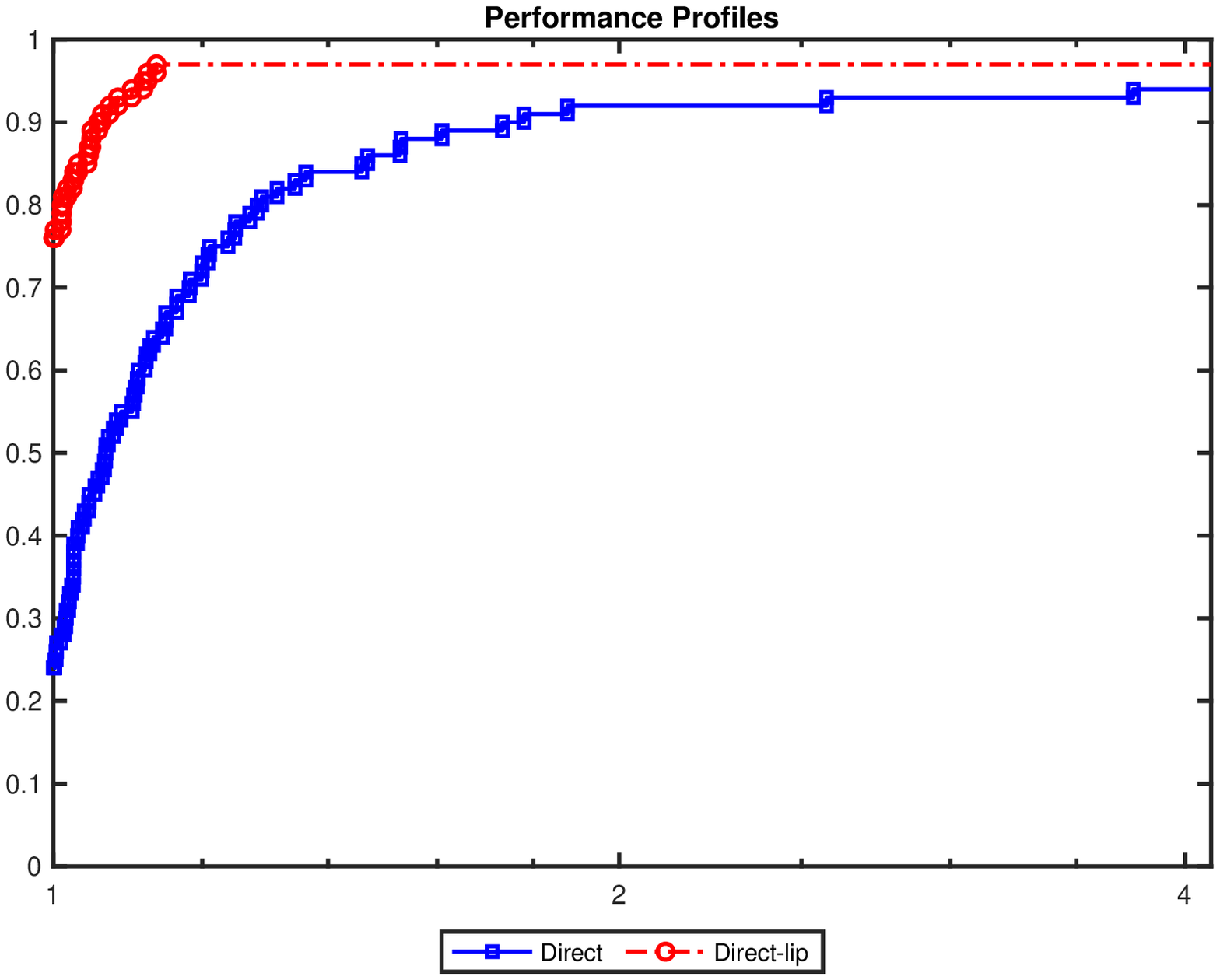}
		\includegraphics[width=0.69\textwidth]{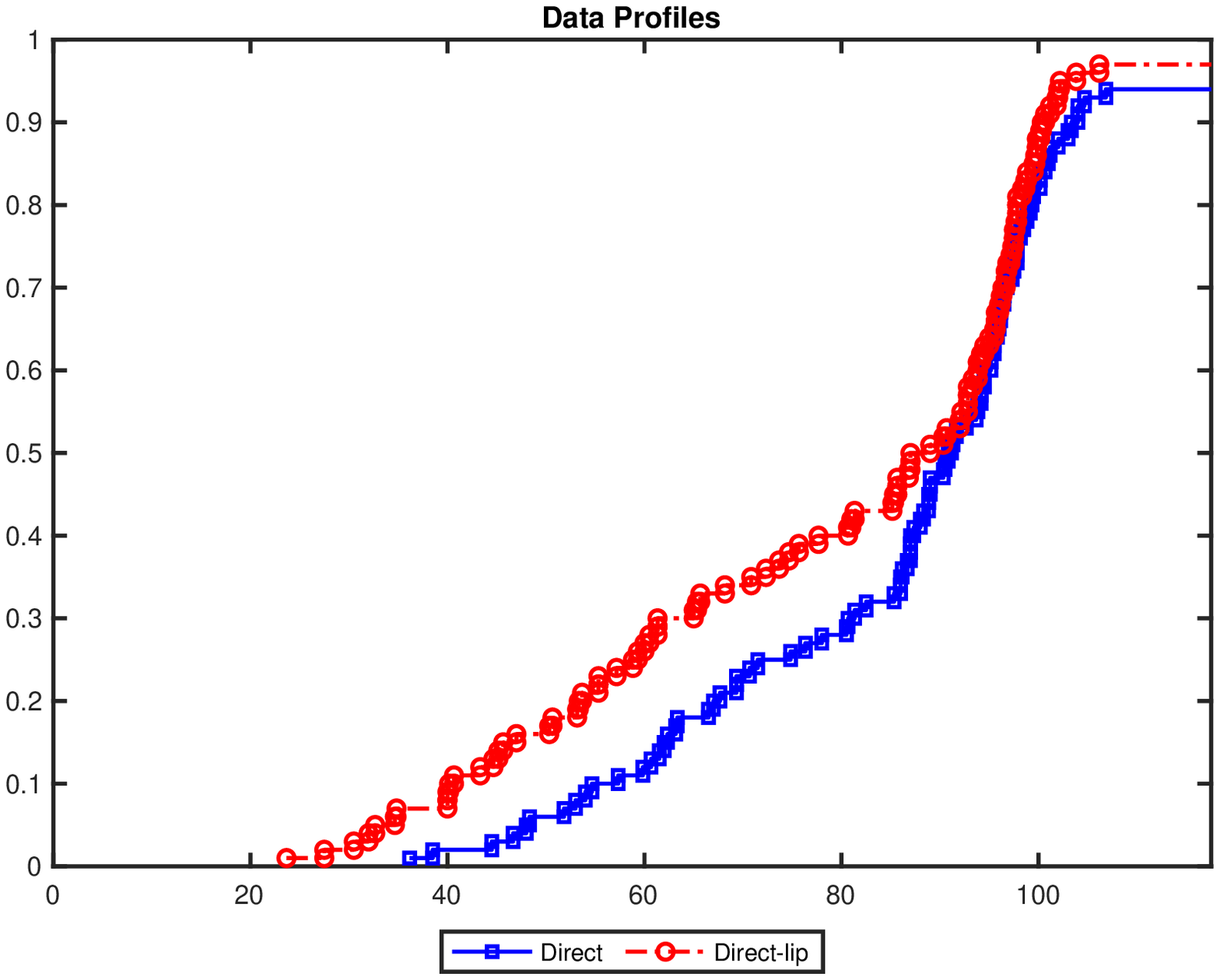}
	\end{center}
	\caption{Performance and data profiles for randomly generated VI problems with trigonometric terms.}
	\label{fig:3}       
\end{figure}

In this subsection we report the results obtained on randomly generated  VI problems with trigonometric terms. We generated 100 instances with 5 variables in this case as well. For each instance, the operator $F(x)=Px+r+T(x)$, where $T_i(x) = w_i \sin(v_i x_i)$, for $i=1,\dots,n$, was randomly built by choosing a matrix $P$ with uniformly distributed random numbers in the interval $[0,3]$ and  vectors $w,v$ and $r$  with uniformly distributed random numbers in the interval $(0,4]$, $(0,2]$, and $[-2,2]$, respectively. 
The box constraints $\{x \in \R^n:\ l\leq x\leq u\}$ were generated by 
considering two vectors $l$ and $u$  with uniformly distributed random
numbers in the interval $[-2,0]$ and $[1,3]$, respectively.
We used the same budget of function evaluation given for affine VIs.
Performance and data profiles are reported in~\autoref{fig:3}. It is easy to see, by taking a look at the  performance profile plot, that the \texttt{$\bar L$-DIRECT} (red line) is again more efficient than \texttt{DIRECT} (blue line), since it gives better performance and satisfies the stopping condition with a smaller number of function evaluations for about the $75\%$ of the instances, and also more reliable. Data profiles show that \texttt{$\bar L$-DIRECT} solves a higher number of instances no matter what the budget used is.

\subsection{Results on VI Problems from the Literature}

\begin{figure}[htbp]
	\begin{center}
		\psfrag{f}{\scalebox{.9}{$\varphi$}}
		\psfrag{min}{\scalebox{.6}{min}}
		\psfrag{Direct}{\scalebox{.6}{\texttt{DIRECT}}}
		\psfrag{Direct-lip}{\scalebox{.6}{$\bar L$-\texttt{DIRECT}}}
		\includegraphics[width=0.69\textwidth]{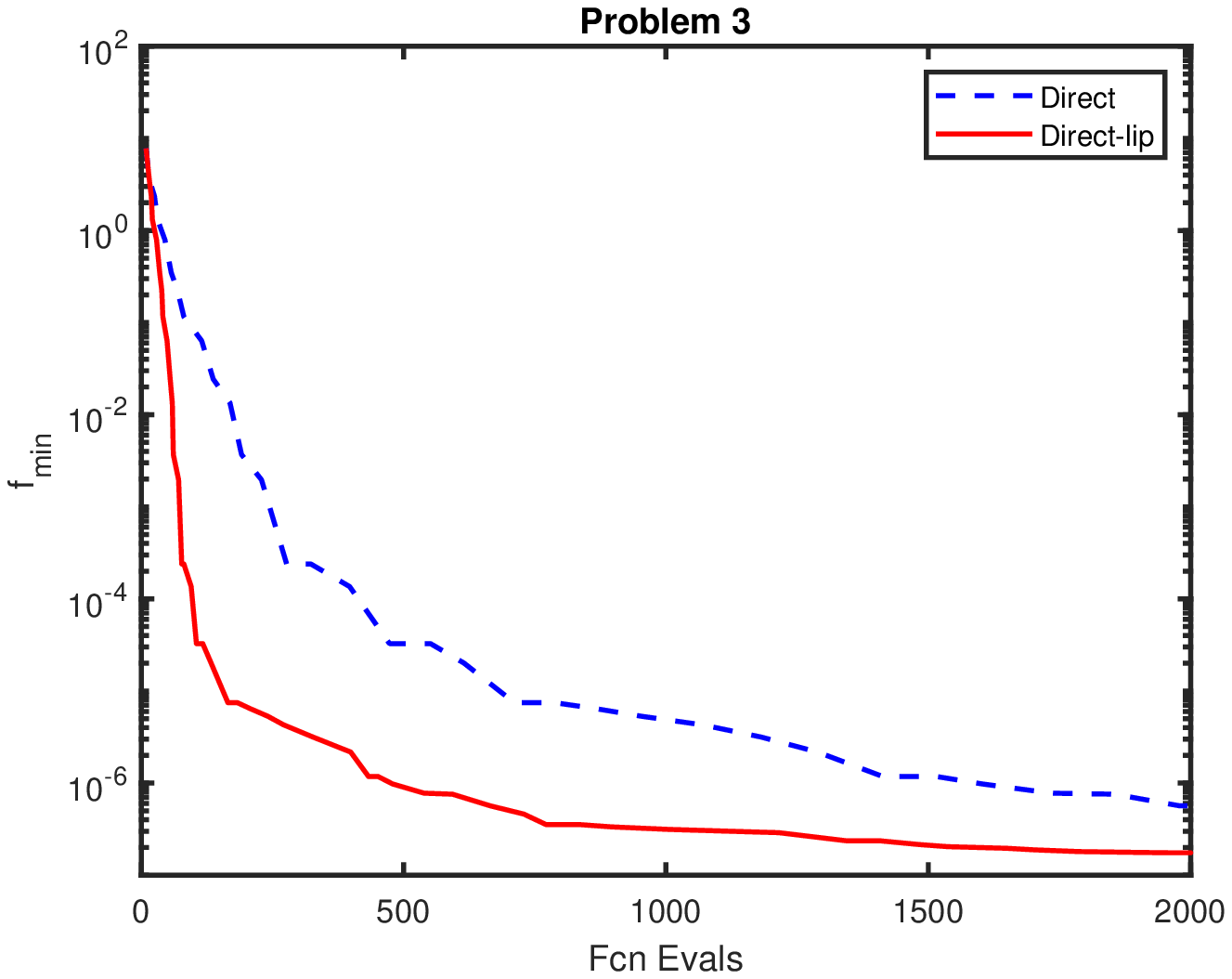}
		\\
		\includegraphics[width=0.69\textwidth]{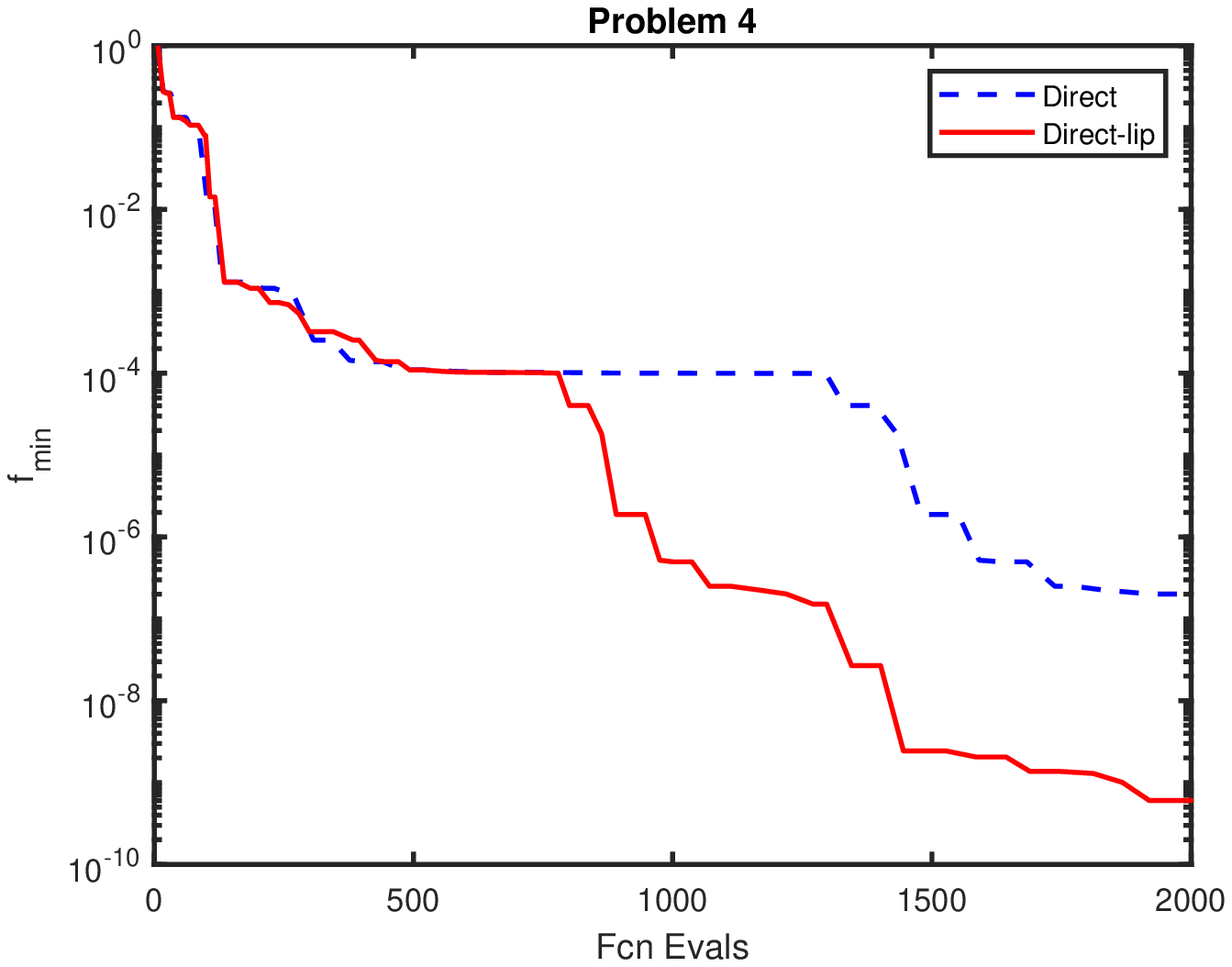}
	\end{center}
	\caption{Comparison between \texttt{DIRECT} and \texttt{$\bar L$-DIRECT} on Problems 3 and 4 from \cite{majig2009global}.}
	\label{fig:1}       
\end{figure}

We finally show results on the Problems 2--6 from paper~\cite{majig2009global}. In order to consider only affine VIs, we dropped the absolute value in the operator $F(x)$ of Problems 4 and 5. 
In \autoref{tab:1}, we report, for each problem, the number of function evaluations needed by the two algorithms to reach a certain gap value (we chose $10^{-1}$, $10^{-3}$, $10^{-5}$).
As we can easily see, the number of function evaluations is usually smaller for \texttt{$\bar L$-DIRECT} (we report in red the cases where \texttt{$\bar L$-DIRECT} needs a higher number of evaluations). In~\autoref{fig:1}, we further report the plots related to the gap reduction with respect to the number of function evaluations used for Problems 3 and 4. We indicate with $\varphi_{min}$ the gap value (reported on the $y$ axis) and with \texttt{Fcn Evals} the number of function evaluations (reported on the $x$ axis). As we can see, the use of the Lipschitz constant estimate significantly speeds up the algorithm.

\begin{table}
	\begin{center}
	\caption{Comparison between \texttt{DIRECT} and \texttt{$\bar L$-DIRECT} on VI Problems from \cite{majig2009global} (\# of f.e. to reach a given gap)}
	\label{tab:1}       
	\begin{tabular}{cc lll lll}
		\hline\noalign{\smallskip}
		
		&  & \multicolumn{3}{c}{\texttt{DIRECT}}&\multicolumn{3}{c}{\texttt{$\bar L$-DIRECT}}   \\
		\noalign{\smallskip}\hline\noalign{\smallskip}
		Problem&$n$ &$10^{-1}$ & $10^{-3}$ &$10^{-5}$ &$10^{-1}$ &$10^{-3}$ &$10^{-5}$\\ 
		\noalign{\smallskip}\hline\noalign{\smallskip}
		\hline
		2 &3&125&431&977&97&313&749\\
		3 &4&115&277&711&49&77&165\\
		4 &3& 79&269&1479&{\color{red}97}& 223&891\\
		5 &5&383&1913&1969&365&{\color{red}1961}&{\color{red}2000}\\
		6 &10&1987&1987&1987&875&1341&1799\\
		\noalign{\smallskip}\hline
	\end{tabular}
	\end{center}
\end{table}


\section{Conclusions}
\label{s:conclusions}

In this paper, we propose a  global optimization approach for solving general EPs without assuming any monotonicity-type condition on $f$. 
This approach is based on two phases: (i) reformulate an EP as a global optimization problem via gap functions; (ii) use an improved version of the \texttt{DIRECT} algorithm, which exploits local bounds of the Lipschitz constant of the objective function, combined with local searches to solve the considered global optimization problem.
Moreover, we provide some general results on Lipschitz continuity of gap functions and, for some special classes of EPs, show simple estimates of their Lipschitz constants that can be exploited in the improved \texttt{DIRECT} algorithm.
Preliminary numerical experiments on a set of instances from the literature and sets of randomly generated instances show the effectiveness of our approach for solving non-monotone EPs.


\begin{thank*}
	The work of M. Passacantando and F. Rinaldi has been partially supported by the Italian Government project PRIN2015B5F27W ``Nonlinear and Combinatorial Aspects of Complex Networks''.
	M. Passacantando is member of the Gruppo Nazionale per l'Analisi Matematica, la Probabilit\`a
	e le loro Applicazioni (GNAMPA - National Group for Mathematical Analysis, Probability and their
	Applications) of the Istituto Nazionale di Alta Matematica (INdAM - National Institute of Higher Mathematics).
\end{thank*}


%
%


\bibliographystyle{plain}
\bibliography{LuPaRi20}   

\end{document}